\documentclass[11pt,reqno]{amsart}

\usepackage{amsmath, amsfonts, amssymb, amscd, enumerate, mathrsfs}
\theoremstyle{plain}
\newtheorem{theo}{Theorem}[section]
\newtheorem{lem}[theo]{Lemma}
\newtheorem{prop}[theo]{Proposition}

\theoremstyle{definition}

\newtheorem{definition}[theo]{Definition}

\newenvironment{pf}{\noindent{\it Proof. }}{$\square$\par\medskip}

\theoremstyle{plain}

\theoremstyle{definition}

\renewcommand{\=}{:=}

\newcommand{\rank}{\operatorname{rank}}
\newcommand{\beq}{\begin{equation}}
\newcommand{\eeq}{\end{equation}}

\renewcommand{\d}{\delta}

\newcommand{\g}{\gamma}
\newcommand{\h}{\eta}

\newcommand{\G}{\Gamma}

\renewcommand{\O}{\Omega}

\newcommand{\bC}{\mathbb{C}}

\newcommand{\bR}{\mathbb{R}}

\newcommand{\bJ}{\mathbb{J}}

\newcommand{\bN}{\mathbb{N}}

\newcommand\U{\mathrm{U}}

\newcommand{\cC}{\mathcal{C}}

\newcommand{\cE}{\mathcal{E}}

\newcommand{\cU}{\mathcal{U}}
\newcommand{\cV}{\mathcal{V}}
\newcommand{\cW}{\mathcal{W}}

\renewcommand{\square}{\kern1pt\vbox
{\hrule height 0.6pt\hbox{\vrule width 0.6pt\hskip 3pt
\vbox{\vskip 6pt}\hskip 3pt\vrule width 0.6pt}\hrule height0.6pt}\kern1pt}

\DeclareMathOperator\Aut{Aut}

\DeclareMathOperator\Id{Id}

\newcommand{\wt}{\widetilde}
\newcommand{\wh}{\widehat}

\newcommand{\be}{\begin{equation}}
\newcommand{\ee}{\end{equation}}

\def\<#1,#2>{\langle\,#1,\,#2\,\rangle}
\newcommand{\arr}{\begin{array}{rlll}}
\newcommand{\ea}{\end{array}}
\newcommand{\bea}{\begin{eqnarray}}
\newcommand{\eea}{\end{eqnarray}}
\newcommand{\bean}{\begin{eqnarray*}}
\newcommand{\eean}{\end{eqnarray*}}

\def\sideremark#1{\ifvmode\leavevmode\fi\vadjust{
\vbox to0pt{\hbox to 0pt{\hskip\hsize\hskip1em
\vbox{\hsize3cm\tiny\raggedright\pretolerance10000
\noindent #1\hfill}\hss}\vbox to8pt{\vfil}\vss}}}

\newcounter{ssig}
\newcounter{ttig}
\newcommand{\refl}{{(\text{\rm ref})}}
\newcommand{\longcorr}{\,\,-\!\!\!\multimap}

   \title[Proper holomorphic maps from pseudoellipsoids]{On  the singular loci  and the images \\ of 
    proper holomorphic maps\\ from pseudoellipsoids}

\begin{document}

  \author{ Cristina Giannotti
and Andrea Spiro }

   \address
{\newline Cristina Giannotti and Andrea Spiro, 
Scuola di Scienze e Tecnologie, Universit\`a di Camerino, Via Madonna delle Carceri 9,
I-62032 Camerino (Macerata),
ITALY\newline
\phantom{a}}

\email
{cristina.giannotti@unicam.it}\par
\medskip

 \email
{andrea.spiro@unicam.it}\par
\medskip


\thanks{{\it Acknowledgments}. This research was partially supported by the Project MIUR ``Real and Complex Manifolds: Geometry, Topology and  Harmonic Analysis'' and by GNAMPA and GNSAGA of INdAM}

\begin{abstract} We prove a generalisation  of    
Rudin's theorem on   proper holomorphic maps 
from the unit ball to the case of proper holomorphic maps 
from pseudoellipsoids.  
\end{abstract}

 \maketitle
\section{Introduction}
In the beginning of the '80's, W. Rudin proved a theorem that gives  an exhaustive  description of proper holomorphic maps $F: B^n \longrightarrow \O$, from the unit ball $B^n$ onto a domain $\O$ of $ \bC^n$, in terms 
of finite unitary reflection groups.  \par
Such result  can be stated as follows.  Recall that for any   finite group $\G$ of automorphisms   of the unit ball there exists  some $h \in \Aut(B_n)$ such that $\G_o = h \G\, h^{-1}$ is a  finite subgroup   of the unitary group $\U_n$,  i.e. of   the group  of  automorphisms  of $B^n$ fixing the origin. Let us denote by 
  $\G_{o\refl} \subset \G_o$   the maximal subgroup of reflections in $\G_o$  and by $(P_1, \ldots, P_n)$ a fixed  set of generators   for the  space of $\G_{o\refl}$-invariant polynomials in $n$-variables.  One can check that  the  holomorphic map 
$$P_{\G}= P_{\G_o} \circ h^{-1}: B^n \longrightarrow B^n_\G \=  P_\G(B^n) \ ,$$
with $P_{\G_o} \= (P_1, \ldots, P_n)$,  is  proper  and   is uniquely associated with  $\G$, up to composition with an element of a special group of polynomial biholomorphisms (see \S 2). 
 Rudin's theorem is the  following.
\par
\begin{theo}[\cite{Ru2}] \label{rudino}ÊFor any proper holomorphic map $F: B^n \longrightarrow \O$  onto a domain $\O \subset \bC^n$, of multiplicity $m > 1$ and $\cC^1$ up to the boundary, there exists a  finite subgroup  $\G \subset\Aut(B^n)$ and a biholomorphism $\Psi: \O \longrightarrow P_\G(B^n)$ such that 
$\Psi \circ F = P_\G $.
\end{theo}
This immediately implies that   any domain, which  is image of  a proper holomorphic map from $B^n$ that is $\cC^1$ on $\overline{B^n}$, is necessarily biholomorphic to  one of the domains $B^n_\G =   P_\G(B^n)$, whose classification can be reduced to  that  of  finite reflection subgroups of $\U_n$.  \par
\smallskip
A crucial element  of Rudin's proof is    the celebrated Alexander Theorem on  global extendability of local automorphisms of $B^n$. One can therefore  ask  if a   result, similar to Rudin's theorem,  can be proved   for  the pseudoellipsoids of $\bC^n$,  on  which   several  analogues of properties of  the unit ball have been obtained   by  appropriate applications of Alexander Theorem (see e.g. \cite{Ru, La, DS, LS}). 
 \par
 \smallskip
So, let us  focus on the pseudoellipsoids of $\bC^n$, namely the domains 
$\cE_{(p)}^n$,  with   $p = (p_1, \ldots, p_k) \in \bN^k$,   $p_i \geq 2$, defined by
 $$\cE^n_{(p)} \= \{\ z \in \bC^n\ :\
 \sum_{j= 1}^{n - k}|z_j|^2 + |z_{n-k+1}|^{2 p_1} + \ldots + |z_n|^{2 p_{k}} < 1\ \}\ .$$
Let also denote by $\varphi^{(p)}: \bC^n \longrightarrow \bC^n$ the holomorphic map 
\beq \label{varphi}   \varphi^{(p)}(z) = (z_1, \ldots, z_{n-k}, (z_{n-k+1})^{p_1}, \ldots, (z_n)^{p_k})\ ,\eeq
whose restriction $\varphi^{(p)}|_{\cE^n_{(p)}}$: $\cE^n_{(p)}$ $\to B^n$ is directly seen to be a proper map. \par
\smallskip
Some  ideas of Rudin's theorem can be actually implemented to  study proper maps from pseudoellipsoids and they bring to
 the following theorem. \par
\begin{theo} \label{budino}ÊFor any proper holomorphic map $F: \cE^n_{(p)} \longrightarrow \O$  onto a domain $\O \subset \bC^n$, of multiplicity $m > 1$ and $\cC^1$ up to the boundary,  there exists a finite subgroup  $\G \subset\Aut(B^n)$   and a proper holomorphic map   $\Psi: \O \longrightarrow P_\G(B^n)$ such that 
$\Psi \circ F = P_\G \circ \varphi^{(p)}$.
\end{theo}
In other words, if we call {\it factoring of  $f$} any expression of the form $f = g \circ h$, where   $f$ appears as composition  of  two  {\it factors}  $g$, $h$,  our theorem  says 
that any  proper holomorphic map $F$,  defined on a pseudoellipsoid and $\cC^1$ up to the boundary, is always a factor of   a map of the form $P_\G \circ \varphi^{(p)}$. This   reduces the analysis of  the first  to that of   factorings  of   the  second.\par
\smallskip
We would like to stress that  our result is optimal, in the sense that  one cannot expect that  $\Psi$ can be proved to be a  biholomorphism, as in Rudin's theorem: 
just consider  the case $F = \Id_{\cE^n_{(p)}}: \cE^n_{(p)}\longrightarrow \cE^n_{(p)}$.\par
 It is also clear that there exist several proper maps   $F$ that are not  equivalent to the trivial examples $ \Id_{\cE^n_{(p)}}$, $ \varphi^{(p)}$ or $P_\G \circ \varphi^{(p)}$. Consider for instance 
 the pseudoellipsoid 
 $$\cE^4_{(2,2)} =   \{\ 
 |z_1|^2 +  |z_2|^2 + |z_3|^4 + |z_4|^4 < 1\ \}$$
  and the map 
 $$F: \cE^4_{(2,2)} \longrightarrow \O = F(\cE^4_{(2,2)})\ ,\qquad F(z) = (z_1 z_2, z_1 + z_2, (z_3)^2, z_4)\ ,$$
 which   is a non trivial factor of the map $P_\G \circ \varphi^{(2,2)}$, given by  
 $$\varphi^{(2,2)}(z) = (z_1, z_2, (z_3)^2, (z_4)^2)\quad \text{and}\quad P_\G(z) = (z_1 z_2, z_1 + z_2, z_3, z_4)$$
 $$ (\text{here}\  P_\G\  \ \text{is associated with the group} \ \G = \{\ \Id_{B^4}\ ,\ g(z) = (z_2, z_1, z_3, z_4)\ \}\ )\ .$$
 \par
 \medskip
 Nonetheless, the fact that $F$ is always a factor of $P_\G \circ \varphi^{(p)}$ gives precise information on the singular locus $Z_F = \{\det J_F(z) = 0  \}$. In fact,  it is necessarily  an analytic subvariety of $\cE^n_{(p)}$ mapped  by $\varphi^{(p)}$  into a subvariety of $B^n$ contained  in  the union  of  the   hyperplanes $\{ z_i = 0\}$ and the fixed point set  of  a finite reflection subgroup of $\Aut(B^n)$.\par
  It also gives strong restrictions on  the class of the images  $\O$ of the proper holomorphic maps from pseudoellipsoids, since,  in their  turn, they are constrained  to admit a  proper holomorphic map onto  a domain  $B^n_\G$. We believe that  such  information can bring to   the classification of  such domains at least  in the most simple cases, as for instance when  $\G$ is trivial and $B^n_\G = B^n$ (see e.g. \cite{BB, KLS, KS}  for the case $n = 2$). \par
\par
\medskip
We finally note that,  when $\cE^n_{(p)} = B^n$,  by Rudin's theorem the  map $\Psi: \O \longrightarrow B^n_\G$,  given in  Theorem \ref{budino},  is necessarily invertible and  the 
 holomorphic correspondence  
$$\Psi^{-1} = F \circ \varphi^{(p)}{}^{-1} \circ P_\G^{-1}: B^n_\G \longcorr \O$$ 
  splits. Therefore a question  worth of further investigations could be  under which   conditions on $\G$ or on  $Z_F$ one can infer that    $\Psi^{-1}$ necessarily splits or, equivalently,  that $\Psi$ is actually a biholomorphism.  \par
 \bigskip
 After a section of preliminaries,   in \S 3 we prove  a crucial property of  the proper holomorphic maps  $F: \cE^n_{(p)} \longrightarrow \O$ that are $\cC^1$ up to the boundary, namely  we show that the  subsets of $B^n$ of the form $\varphi^{(p)}(F^{-1}(w))$, $w \in \O$,  coincide with  the orbits  of a  finite group $\G$ of automorphisms of $B^n$.  With the help of this fact, we  prove   Theorem \ref{budino} in \S 4.\par
\medskip
\noindent {\bf Acknowledgments.} We are  grateful to the referee for  his/her  kind and  helpful remarks.
\section{Preliminaries}
 \subsection{Finite subgroups  $\G \subset \Aut(B^n)$ and the  proper maps $P_\G$}\hfill\par
 \label{section2.5}
 Let us  call {\it geodesic hyperplane of $B^n$} any $(n-1)$-dimensional    subvariety of  $B^n$ of  the form  
$g\left( \{\ z \in B^n\ : \ z_n = 0\ \}\right)$ for some  $g \in \Aut(B^n)$.
Note that the  geodesic  hyperplanes $g(\{z_n = 0\})$, determined by elements $g \in \U_n = \Aut_0(B^n)$,  are  the usual affine hyperplanes through the origin. 
\par
  \begin{definition} 
 An element $h \in \Aut(B^n)$ is called Ê{\it  reflection} if  it has finite period and its fixed point set  is a  geodesic hyperplane.  
 \end{definition}
For a given finite subgroup $\G \subset \Aut(B^n)$, we denote by $\G_{\refl} \subset \G$ the   subgroup  generated by  all reflections in $\G$. Note that 
a non trivial element $g \in \G$  is a reflection if and only if  its fixed point set is  $(n-1)$-dimensional  (in fact, 
 up to conjugation in $\Aut(B^n)$, any such element is in $\U_n$). This implies that $\G_\refl$ is normal in  $\G$. \par
 \bigskip
Consider a finite reflection group  $\G_o$ $= \G_{o\refl}$ in  $ \U_n$. 
By a classical result of  Chevalley  (\cite{Ch, ST, Fl}), there are $n$ homogeneous, $\G_{o}$-invariant polynomials  $P_1$, \ldots, $P_n$  that constitute   a basis
for the invariants of   $\G_{o}$ (i.e,  the  $\G_{o}$-invariant polynomials   $f \in \bC[z_1, \ldots, z_n]$  are exactly those  of the form  $f = q(P_1, \ldots, P_n)$ for some  $q\in \bC[z_1, \ldots, z_n]$).  The map  
$$P_{\G_{o}}(z) = (P_1, \ldots, P_n): B^n \longrightarrow \bC^n$$
is  uniquely determined by  $\G_{o}$,   up to composition with  the polynomial maps that interchange the bases of homogeneous polynomials  for the invariants $\G_{o}$. The group of such basis changes   is the same for all  groups  $\G_{o}'$ of  the  conjugacy class of  $\G_o = \G_{o\refl}$ in $\U_n$.  \par
\bigskip
Consider now   an arbitrary finite group of automorphisms $\G\subset \Aut(B^n)$, with  reflections  subgroup $\G_\refl$. 
 It is known that the elements of $\G$ have a common fixed point $x_o$ (see e.g. \cite{Ru2}, Thm. 3.1), so that for any  $h \in \Aut(B^n)$ with   $h(x_o) = 0$, the conjugate group 
 $\G_o = h  \G  h^{-1} $ is in $\U_n$ and has 
  $ \G_{o\refl} = h  \G_\refl  h^{-1}$ as   reflection subgroup. 
 We may therefore consider the   map 
$$P_\G: B^n \longrightarrow \bC^n\ ,\qquad P_\G = P_{\G_{o\refl}} \circ h\ ,$$
whose components are $\G_\refl$-invariant rational functions. Up to compositions with the basis changes described above,  $P_\G$ is  uniquely determined by  $\G_\refl$.
By \cite{Ru2}, Thm. 2.5, 
the image $B^n_\G = P_\G(B^n)$
is   a  domain of $\bC^n$, which   is   uniquely determined by the subgroup $\G_\refl\subset \G$  up to  biholomorphisms, and   $P_\G: B^n \longrightarrow B^n_\G$  is  a proper holomorphic map.  \par
\bigskip
We conclude recalling    the statement of Rudin's generalisation of Alexander Theorem. Let us  call   {\it local automorphism of $B^n$\/}  any
 biholomorphism $f: \cU_1 \subset B^n \longrightarrow \cU_2 \subset B^n$ between  connected open subsets of $B^n$ such that:
 \begin{itemize}
 \item[a)] each of the intersections $\partial  \cU_i \cap \partial  B^n$, $i = 1,2$, contains a boundary  open set
   $\Gamma_i \subset \partial B^n$;
  \item[b)] there exists a sequence $\{x_k\} \subset \cU_1$ which converges to a point $x_o \in \Gamma_1$, which is  not a limit point of  $\partial \cU_1 \cap  B^n$, and so that $\{f(x_k)\}$
  converges  to a point $\hat x_o \in \Gamma_2$, which is  not a limit point of  $\partial \cU_2 \cap  B^n$.
 \end{itemize}
Let also  say that  $f$  {\it extends   to a global automorphism} if there exists
 $F \in \Aut(B^n)$ such that $F|_{\cU_1} = f$. \par
 \begin{theo}[\cite{Al,Ru}] \label{alexander}ÊAny local automorphism of $B^n$ extends to a global one.
 \end{theo}
\bigskip
 \subsection{Correspondences}\hfill\par
 Let $D, D' \subset \bC^n$ be two bounded domains and denote by
 $\pi: D \times D' \longrightarrow D$ and  $\pi': D \times D' \longrightarrow D'$
 the two natural projections.  We recall that a  {\it holomorphic correspondence between $D$ and $D'$} is a
 subvariety $V \subset D \times D'$. It is called  {\it proper} if  the  restricted projections   $\pi|_V: V \longrightarrow D$ and $\pi'|_V: V \longrightarrow D'$
 are proper maps.  A holomorphic correspondence is called {\it irreducible}
 if it is an irreducible subvariety.  \par
 \smallskip
  A  holomorphic  correspondence $V$  is uniquely determined by the associated multivalued map 
 $$f: D \longcorr D'\ ,\qquad f(z) = \pi'\left(\pi|_V^{-1}(z)\right)\ ,$$
which  is a  (single-valued) holomorphic map   if and only if   $\pi|_V$ is   injective.
We  often denote a  holomorphic correspondence $V$ by  the corresponding  multivalued map $f$, so that the subvariety $V$ coincides with the graph
$$V  = \G_f := \{\ (z, w)\in D \times D'\ :\ w \in f(z)\ \}\ .$$
If $f: D \longcorr D'$ is a  holomorphic correspondence, we denote by $f^{-1}: D' \longcorr D$ the holomorphic correspondence with
$$\G_{f^{-1}} =  \{\ (w, z)\in D' \times D\ :\  (z, w) \in \G_f\ \}\ .$$
If $f: D \longcorr D'$ and $f': D' \longcorr D''$ are two (proper) holomorphic correspondences, it is known that 
the multivalued map $f' \circ f: D \longcorr D''$ with 
$$\G_{f' \circ f} =  \{\ (z, v)\in D \times D''\ :\  (z, w) \in \G_f\ ,\ (w, v) \in \G_{f'}\ \text{for some}\ w \in D'\  \}$$
is a (proper) holomorphic correspondence as well (\cite{St}). \par
Finally, given two (proper) holomorphic correspondences $f_1$, $f_2: D \longcorr D'$, we denote by 
$f = f_1 \cup f_2$
 the (proper) holomorphic correspondence  with  $\G_f = \G_{f_1} \cup \G_{f_2} \subset D \times D'$.  \par
 \smallskip
 If $f: D \longcorr D'$ is a proper holomorphic correspondence, there exist a positive integer $p$
 and a subvariety $W \subset D$ such that, for every $z_o \in D \setminus W$, there are an open neighbourhood $\cU \subset D \setminus W$ of $z_o$ 
 and $p$ holomorphic maps $f_i: \cU \longrightarrow D'$ such that the sets $f(z)$, $z \in \cU$,  have cardinality $p$ and are equal to
 $$f(z) = \{\ f_1(z),\  \ldots, \ f_p(z)\ \}\ .$$
We shortly say that ``$f$ is a $p$-valued map'' . 
For a given $z_o \in D$, we say that  {\it $f$ splits at $z_o$}  if there exists a neighbourhood $\cU \subset \bC^n$ of $z_o$ such that 
$$\G_{f|_{\cU}} = \G_f \cap \pi^{-1}(\cU \cap D) = \G_{f_1} \cup \ldots \cup \G_{f_q}$$
for  some  single-valued holomorphic maps  $f_i: D \cap \cU \longrightarrow D'$. If $f$ is $p$-valued, the number of such  single-valued maps has to coincide with   $p$.\par
\smallskip
  We say that {\it $f$ splits} if it splits at all points. If $D$ is simply connected, $f$ splits  if and only if 
  there are $p$  holomorphic maps $f_i: D \longrightarrow D'$ such that 
   $f = f_1  \cup \ldots \cup f_p$. The  $f_i$'s  are called  {\it single-valued components of $f$}.\par
   The following  is a direct consequence of  \cite{BB}, Lemma 3.1. 
   \begin{lem} \label{splittinglemma}ÊIf   $f: D \longcorr D'$ is a holomorphic correspondence, either  it splits or 
   there exists an analytic subvariety $S_f \subset D$ of  dimension $n-1$ such that $f$ does not split at $z$  for any  $z \in S_f$.
   \end{lem}
\bigskip
 
\subsection{A technical fact concerning proper holomorphic maps }\hfill\par
Let   $F: D \longrightarrow D'$ be a  proper holomorphic map  with multiplicity $m$ and denote $Z_F = \{x \in D\ : \det J_F = 0\}$.  
  If  $F$  extends to a $\cC^1$-map  $F: \cU  \longrightarrow \bC^n = \bR^{2n}$  on a neighbourhood $\cU$ of $\overline D$, we denote by $\bJ_F(x)$, $x \in \cU$,  the (real) Jacobian of $F$ at $x$, where $F$ is considered as a map between open subsets of $\bR^{2n}$. 
If  such (real)  map is expressed in terms of the complex coordinates $(z_i, \overline z_i)$  and $F$ is holomorphic at $x$, then
 \beq \label{lillina}Ê\bJ_F(x) = \left(\begin{array}{cc}  \frac{\partial F_i}{\partial z_j}  & 0 \\ 0 & \overline{ \frac{\partial F_i}{\partial z_j} }\end{array}\right)  \quad \text{and hence} \quad \rank \bJ_F(x) = 2 \rank \left( \frac{\partial F_i}{\partial z_j} \right) \ .\eeq
By continuity,  \eqref{lillina}  holds in  $\overline D$, so  that  $\rank (\bJ_F(x))$ is even  for all $x \in \overline D$.\par
\smallskip
\begin{lem}\label{lilla1}  Let   $D \subset \bC^n$ be a  bounded domain  with smooth boundary and $F: D \longrightarrow D'$ a proper holomorphic map  admitting a $\cC^1$ extension to  $\overline D$. Let us also use the notation   $Z_{F} \cap \partial D\= \{x \in \partial D \ : \det \bJ_F = 0\}$. \par
Then, the $(2n-1)$-dimensional Hausdorff measure of $F(Z_F\cap \partial D)$ is $0$.
\end{lem}
\begin{pf} Let $x \in  \partial D$ and consider a system of real coordinates $\xi = (x_1, \ldots, x_{2n})$ on a neighbourhood $\cV$ of  $x$ such that $\partial D \cap \cV = \{Ê\ x_{2n} = 0\ \}$. In such coordinates, the restriction $\wt F = F|_{\partial D}$ is of the form 
$\wt F(x_1, \ldots, x_{2n-1}) = F(x_1, \ldots, x_{2n-1}, 0)$ and the Jacobian  $\bJ_F(x)$ of $F$ is of the form
$$\bJ_F(x) = \left(\begin{matrix} \ \ \ J_{ \wt F}(x) \ \ \ & \left|\begin{matrix} \frac{\partial F_1}{\partial x_{2n}} (x)\\
\vdots
\\
\frac{\partial F_{2n}}{\partial x_{2n}} (x)
 \end{matrix}\right.\end{matrix}\right)\ .$$
This means that  $\rank \wt F|_x \leq \rank F|_x \leq \rank \wt F|_x + 1$.  If $x \in Z_F \cap \partial D$,  previous remarks imply that  $ \rank \wt F|_x \leq \rank F|_x\leq 2n-2$ and, conversely, if  $\rank \wt F|_x  \leq 2n-2$, one has that  $\rank F|_x  \leq 2n-1$ and hence  $x \in Z_F\cap \partial D$.  This means that  $Z_{F} \cap \partial DÊ= \{\ x \in \partial D\ : \ \rank \wt F|_x \leq 2n-2\ \}$ and 
the  claim follows from generalised  Morse-Sard Theorem (see e.g. \cite{M}). 
\end{pf} 
\par
 \bigskip
\section{$F$-related points in $\cE^n_{(p)}$ and $B^n$}
In all the following,   $F: \cE_{(p)}^n \longrightarrow \O \subset \bC^n$ is a  proper holomorphic map of multiplicity $m$ and 
$\varphi^{(p)}: \cE_{(p)}^n \longrightarrow B^n$ is the proper holomorphic  map defined in \eqref{varphi}.  We also set 
\beq \label{pi} \pi = \{\ z \in \bC^n\ :\ z_{n - k + 1} \cdot z_{n-k + 2} \cdot \ldots \cdot z_n = 0\ \}\ .\eeq
\par
\begin{definition} A  subset $J \subset \cE^n_{(p)}$  is called {\it complete $F$-set in $\cE^n_{(p)}$} if
$J  = F^{-1}(w_o)$
for some $w_o \in \O$.  It  is called {\it good} if it is the pre-image of a point
$$w_o \in \O \setminus F\left(Z_F \cup \pi\right)\ .$$
Similarly, a  subset $\wt J \subset B^n$  is called {\it  complete $F$-set in $B^n$} if it is of the form
$\wt J  =  \varphi^{(p)}(J)$
 for a complete $F$-set in $\cE^n_{(p)}$. If $J$ is good,    also  $\wt J$ is called {\it good}. \par
Two points  of a complete $F$-set in $\cE^n_{(p)}$  (resp.  in $B^n$) are called   {\it $F$-related}. Similarly,  two sequences $\{x_k\}$, $\{x'_k\}$ in $ \cE^n_{(p)}$  (resp.  in $B^n$)  are  called  {\it $F$-related} if $x_k$ and $x'_k$  are $F$-related for all $k$'s.
\end{definition}

  \begin{lem} \label{successioni}ÊLet $F: \cE_{(p)}^n \longrightarrow \O \subset \bC^n$ be a proper holomorphic map of multiplicity $m > 1$,  admitting a $\cC^1$ extension to $\overline{\cE^n_{(p)}}$. Then there  exist $m$ pairwise $F$-related sequences $\{x^{(1)}_{k}\} $, \ldots,   $\{x^{(m)}_{k}\}$  in $\cE^n_{(p)}$  with the following properties: 
  \begin{itemize}
  \item[i)] they converge to  $m$ distinct  points $x^{(1)}_o, \ldots, x^{(m)}_{o}  \in \partial \cE^n_{(p)}$; 
  \item[ii)] there are  disjoint connected open sets $\cU^{(i)} \subset \bC_n$,  $1 \leq i \leq m$, such that:
  \begin{itemize}
  \item[--]  $x^{(i)}_{o} \in \cU^{(i)}$; 
  \item[--]Êthe restrictions $F|_{\cU^{(i)}Ê\cap \cE^n_{(p)}}: \cU^{(i)}Ê\cap \cE^n_{(p)}\longrightarrow F(\cU^{(i)}Ê\cap \cE^n_{(p)})$ are biholomorphisms onto the same 
   open set $\cW\! = F(\cU^{(i)}Ê\cap \cE^n_{(p)}) \subset \O$; 
  \item[--] the restrictions  
  $\varphi^{(p)}|_{\cU^{(i)}Ê\cap  \cE^n_{(p)}}: \cU^{(i)}Ê\cap \cE^n_{(p)}\longrightarrow \cV^{(i)}Ê= \varphi^{(p)}( \cU^{(i)}Ê\cap \cE^n_{(p)})$
  are biholomorphisms. 
  \end{itemize}
    \end{itemize}
      \end{lem}
 \begin{pf} Let $\wt Z_F = Z_F \cap \partial \cE^n_{(p)}$ and  
 $\wt \pi = \pi \cap \partial \cE^n_{(p)}$.
 Since $F$ is Lipschitz in $\overline{\cE^n_{(p)}}$ and the $(2n-1)$-dimensional Hausdorff measure $H_{2n-1}(\wt \pi)$ is zero, we have $H_{2n-1}(F(\wt \pi)) = 0$. Hence,  by  Lemma \ref{lilla1},  
  $$H_{2n-1}(F(\wt Z_F \cup \wt \pi)) \leq H_{2n-1}(F(\wt Z_F)) + H_{2n-1} (F(\wt \pi )) = 0 \ .$$
Since $\partial \O$ surely includes pieces of smooth hypersurfaces, $H_{2n-1}(\partial \O) >0$ and  consequently 
$$\partial \O \setminus F(\wt Z_F \cup \wt \pi)  \neq \emptyset\ .$$ 
 Pick a point  $w_o \in \partial \O \setminus F(\wt Z_F  \cup \wt \pi)$,   a pre-image $x^{(1)}_o \in F^{-1}(w_o)$  and a small  arc $\g^{(1)}_t \subset \cE^n_{(p)} \setminus F^{-1}(F(Z_F))$, $t \in [0, 1)$,  ending at $x^{(1)}_o = \lim_{t \to 1} \g^{(1)}_t$. Since the restriction of $F$ to $\cE_{(p)}^n \setminus F^{-1}(F(Z_F))$ is a proper, unbranched  cover  (see e.g. \cite{Be}),   there are exactly  $m-1$ disjoint arcs $\g^{(2)}_t$, \ldots, 
 $\g^{(m)}_t$,  $t \in [0, 1)$,  determined by  the points that are  $F$-related  to the points $\g^{(1)}_t$.\par 
 \smallskip
 By construction,    $x^{(i)}_o = \lim_{t \to 1} \g^{(i)}_t \in F^{-1}(w_o)$ for all $2 \leq i \leq m$. 
 Since $\det \bJ_F(x_o^{(i)}) \neq 0$, any such point admits a connected neighbourhood,  on which $F$ is an homeomorphism, implying  that $x^{(1)}_o$, \ldots, $x^{(m)}_o$ are  all distinct. \par
 We may consider disjoint  connected neighbourhoods $\wh \cU^{(i)}$ of the $x_o^{(i)}$ that do not intersect $Z_F \cup \pi$, so that $F|_{\wh \cU^{(i)}}$ and $\varphi^{(p)}|_{\wh \cU^{(i)}}$
 are homeomorphisms onto their images in $\bC^n$ and the restrictions
  $F|_{\wh \cU^{(i)} \cap \cE^n_{(p)}}$ and  $\varphi^{(p)}|_{\wh \cU^{(i)} \cap \cE^n_{(p)}}$ are biholomorphisms onto their images in $\O$ and $B^n$, respectively. Setting  $\wh \cW = \bigcap_{i = 1}^m F(\wh \cU^{(i)})$ and  $\cU^{(i)}=  F^{-1}(\wh \cW) \cap \wh \cU^{(i)}$, any choice of
  $F$-related  sequences  $\{x^{(i)}_{k} = \g^{(i)}_{t_k}\} $
 with $\lim_{k \to \infty} t_k = 1$  satisfies the  claim. \end{pf}
 \medskip
From now on, we consider a fixed choice of   $m$ sequences $\{x^{(j)}_{k}\} $, $1 \leq j \leq m$, converging to 
$ x^{(j)}_o \in \partial \cE^n_{(p)}$,  and open sets 
\beq \label{cristina}Ê \cU^{(i)}Ê\qquad \text{and}\qquad \cV^{(i)} = \varphi^{(p)}( \cU^{(i)} \cap \cE^n_{(p)})  \subset B^n\eeq
satisfying the statement of  Lemma 
\ref{successioni}.  We also denote by $g^{(i,j)}: \cV^{ (i)} \to  \cV^{ (j)}$ the biholomorphisms
\beq \label{fabrizio} g^{(i, j)} =  \left(\varphi^{(p)}|_{\cU^{(j)} \cap \cE^n_{(p)}}\right) \circ \left( F|_{\cU^{(j)} \cap \cE^n_{(p)}}\right)^{-1} \circ 
 \left( F|_{\cU^{(i)} \cap \cE^n_{(p)}}\right)\circ \left(\varphi^{(p)}|_{\cU^{(i)} \cap \cE^n_{(p)}}\right)^{-1}\ .\eeq
 Notice that {\it the   $g^{(i,j)}$'s are local automorphisms of $B^n$ and, by Theorem  \ref{alexander},   they all  extend  to  global automorphisms of $B^n$}. 
  We finally set 
\beq \label{definitionofgamma}Ê\G = \{\ g^{(i,j)}\ ,\ 1 \leq i, j\leq m\ \} \subset \Aut(B^n)\ .\eeq\par
\medskip
\begin{prop} \label{proposition3.3} Two points $y, y' \in B^n$ are $F$-related in $B^n$ if and only if $y' = g(y)$ for some $g \in \G$. In particular, $\G$ is a finite subgroup of $\Aut(B^n)$.
\end{prop}
\begin{pf}   We first prove the necessity.  Let  $y$, $y' \in B^n$ be $F$-related points, i.e.   $y =  \varphi^{(p)}(x)$,  $y' = \varphi^{(p)}(x')$ for two points
$x$, $x'$ of  a  complete $F$-set $J = \{x_1, \ldots, x_m\} = F^{-1}(w)$  in $\cE^n_{(p)}$.
Let also $Z_{F, \varphi^{(p)}} $ be the analytic subvariety of $\cE^n_{(p)}$ defined by 
$$Z_{F, \varphi^{(p)}} = F^{-1}\left(F\left(Z_F \cup \pi\right)\right) $$
where, as usual,  $\pi = \{z_{n - k + 1} \cdot \ldots \cdot z_n = 0\}$.
We consider two cases.\par
\medskip
\noindent {\it Case 1: $J$ is good, i.e. $J \cap Z_{F, \varphi^{(p)}} = \emptyset$ }. \par
Since $Z_{F, \varphi^{(p)}} $ is  analytic subvariety of $\cE^n_{(p)}$, the set 
$\cE^n_{(p)} \setminus  Z_{F, \varphi^{(p)}}$  is connected (see e.g. \cite{Na}, Ch.4, Prop. 1).
%
%
We may therefore 
consider a $\cC^0$ curve  
  $\h: [0,1] \longrightarrow \overline{\cE^n_{(p)}}$ 
such that 
 \begin{itemize}
 \item[--] $\h_0 = x$ and $\h_1 = x_o^{(1)} $; 
 \item[--]  $\h_t \in  \cE^n_{(p)}\setminus Z_{F, \varphi^{(p)}}$ for any  $0 \leq t < 1$. 
  \end{itemize}
The corresponding curve $\g = \varphi^{(p)} \circ \h: [0,1]\longrightarrow \overline{B^n}$ is such that 
 \begin{itemize}
 \item[--] $\g_0 = y$ and $\g_1 = y_o^{(1)} = \varphi^{(p)}(x_o^{(1)})$; 
 \item[--]   $\g_t \in  B^n\setminus\varphi^{(p)}\left(Z_{F, \varphi^{(p)}}\right)$ for any $0 \leq t<1$. 
   \end{itemize} 
  \par
  \medskip
Consider now a $\cC^0$-curve $\h': [0,1]\longrightarrow \overline{\cE^n_{(p)}}$ such that $Ê\h'_0 = x'$ and 
$\h'_t$ is $F$-related to $\h_t$ for any $0 \leq t < 1$. By the properties of proper holomorphic maps and the fact that $\h_t \notin F^{-1}(F(Z_F))$, such curve exists and it is unique. 
In particular, $\h'_t \in \cE^n_{(p)}\setminus Z_{F, \varphi^{(p)}}$ for any $t < 1$ and $ \h'_1\in \partial \cE^n_{(p)}$. 
  \smallskip
  Finally,  let $\g': [0,1] \longrightarrow \overline{B^n}$ be the curve $\g' = \varphi^{(p)} \circ \h'$. By construction, 
 $$\g'_0 = y'\ ,\qquad \g'_t \in B^n\setminus\varphi^{(p)}\left(Z_{F, \varphi^{(p)}} \right)\  \text{if} \  t< 1\ , \qquad   \g'_1 =  \varphi^{(p)}(\h'_1)\in \partial B^n\ .$$  
 \smallskip
Notice that, being  $\h_t$ and $\h'_t$  distinct and $F$-related, the end-point  $\h'_1$ must be  one of the points $x_o^{(2)}$,  $x_o^{(3)}, \ldots,$ $ x_o^{(m)}$.  For simplicity, we assume  $\h'_1 = x_o^{(2)}$. \par
\medskip
Now, we observe that,  for any $t \in [0,1)$, there exist neighbourhoods $\cU_t$, $\cU'_t \subset \cE^n_{(p)}$ of $\h_t$ and $\h'_t$, respectively, and neighbourhoods $\cV_t$, $\cV_t' \subset B^n$ of $\g_t$ and $\g'_t$, such that the restrictions 
$$
\begin{array}{ll}
F|_{\cU_t}: \cU_t \longrightarrow F(\cU_t)\ ,&\varphi^{(p)}|_{\cU_t}: \cU_t \longrightarrow \cV_t\ ,\\
\ & \\
F|_{\cU'_t}: \cU'_t \longrightarrow F(\cU'_t)  = F(\cU_t)\ ,& \varphi^{(p)}|_{\cU'_t}: \cU'_t \longrightarrow \cV'_t
\end{array}$$
are biholomorphisms, so that  also 
\beq \label{accacont}Êh_t = \varphi^{(p)}|_{\cU'_t}Ê\circ \left(F|_{\cU'_t}\right)^{-1} \circ F|_{\cU_t} \circ \left(\varphi^{(p)}|_{\cU_t}\right)^{-1}: \cV_t \longrightarrow \cV'_t\eeq
is a biholomorphism. 
For    $t = 1$, we  set  $\cV_1=\cV^{(1)}$,  $\cV'_1=\cV^{(2)}$ and 
\beq h_1 =  g^{(1,2)}|_{\cV_1}: \cV_1 \longrightarrow \cV'_1\ .\eeq
\medskip
We claim that, for any $t, s \in [0,1]$,  with $\cV_{t} \cap \cV_{s} \neq \emptyset$,  
\beq \label{3.9} h_t|_{\cV_{t} \cap \cV_{s}}Ê= h_s|_{\cV_{t} \cap \cV_{s}}\ .Ê\eeq
In fact, if $\cV_{t} \cap \cV_{s} \neq \emptyset$ (hence, it contains a subarc of $\g$), then   $\cU_t \cap \cU_{s} \neq \emptyset$  (it contains a subarc of $\h$) and $\varphi^{(p)}|_{\cU_t \cap \cU_{s}}$ is a biholomorphism onto $ \varphi^{(p)}(\cU_t \cap \cU_{s}) $  with  inverse 
$$\left(\varphi^{(p)}|_{\cU_t \cap \cU_{s}}\right)^{-1}Ê= \left.\left(\varphi^{(p)}|_{\cU_t}\right)^{-1}\right|_{\varphi^{(p)}(\cU_t \cap \cU_{s})}Ê= \left.\left(\varphi^{(p)}|_{\cU_s}\right)^{-1}\right|_{\varphi^{(p)}(\cU_t \cap \cU_{s})}\ .$$
By a similar argument 
$$\left(F|_{\cU'_t \cap \cU'_s}\right)^{-1}  = \left.\left(F|_{\cU'_t}\right)^{-1}\right|_{F(\cU'_{t} \cap \cU'_{s})}Ê= \left.\left( F|_{\cU'_s}\right)^{-1}\right|_{F(\cU'_{t} \cap \cU'_{s})}$$
and \eqref{3.9}Ê follows directly from  the definitions of the $h_t$'s. 
\par
\smallskip
By compactness,  there are $t_1$, \ldots, $t_{N-1}$, $t_N=1$ $\in [0,1]$ such that $\g([0,1]) \subset \bigcup_{k = 1}^N \cV_{t_k}$ and, by \eqref{3.9},   the maps $h_{t_i}$ can be  glued together to determine a  holomorphic map 
$$h: \cV =\bigcup_{k = 1}^N \cV_{t_k} \longrightarrow \cV' = \bigcup_{k = 1}^N \cV'_{t_k}\ .$$
Since $h|_{\cV_{1}Ê}= h_1 =g^{(1,2)}|_{\cV_{1}Ê}$, by the Identity Principle, $h = g^{(1,2)}|_{\cV}$ and 
$y' = h(y) = g^{(1,2)}(y)$,  proving the claim.\par
\bigskip
\noindent{\it Case 2: $J$ is not good,  i.e. $J \cap Z_{F, \varphi^{(p)}} \neq \emptyset$}.
In this case  $J = F^{-1}(w)$ for some  $w\in F\left(Z_F \cup \pi \right)$. Let $\{w_k\} \subset
\O \setminus F\left(Z_F \cup\pi \right)$ be a sequence  
 with $\lim_{k \to \infty}Êw_k = w$ and denote by  $\wt J_k = \varphi^{(p)}(F^{-1}(w_k)) = \{y_{k,1}, \ldots, y_{k,r_k}\}$ the corresponding sequence 
of good complete $F$-sets in $B^n$. 
Taking a suitable subsequence, we may assume that $y$,  $y'$ are limits of two sequences $\{y_k\}$, $\{y'_k\}$  with $y_k, y'_k \in \wt J_k$ for any $k$. By the previous part of the proof, there are $g_k \in \G$ such that $g_k(y_k) = y'_k$. 
Since $\G$ is a finite set, we may consider a subsequence $\{y_{k_n}\}$ and $g \in \G$ such that 
$g(y_{k_n}) = y'_{k_n}$ for any $n$. Therefore  $g(y) = \lim_{n \to \infty} g(y_{n_k}) = \lim_{n \to \infty} y'_{n_k} = y'$ and the claim follows.\par
\bigskip
Let us now prove the sufficiency. Let $y, y' \in B^n$ be such that $y' = g(y)$ for some $g \in \G$. If $y = y'$, there is nothing to prove. Therefore, we assume  $y \neq y'$ and  $g \neq \Id_{B^n}$.  For simplicity, we assume that $g = g^{(1,2)}$.  Consider  the analytic subvariety of $\overline{B^n}$
\beq \label{zetaprimo}ÊZ'Ê= \bigcup_{h \in \G} h\left(\varphi^{(p)}(Z_{F, \varphi^{(p)}})\right) \eeq
and we prove the claim in the mutually  exclusive cases $y, y' \notin Z'$  and $y, y' \in Z'$, respectively. 
\par
\bigskip
{\it Case 1: $y, y' \in B^n \setminus Z'$}.   Pick a point $y_o \in \cV^{(1)} \setminus Z' \subset B^n$ and observe that, 
being  $ B^n \setminus Z'$ complementary to an analytic subvariety, there exists a $\cC^0$ curve  $\g: [0,1]\longrightarrow B^n$  such that 
 \begin{itemize}
 \item[--] $\g_0 = y$ and $\g_1 = y_o$; 
 \item[--]   $\g_t \in  B^n\setminus Z'$ for any $t \in [0,1]$; 
  \end{itemize}
  Secondly, consider the $\cC^0$ curve 
 $\g' = g^{(1,2)} \circ \g: [0,1]\longrightarrow B^n$. 
By construction,  
$\g'_0 = y'$ and $\g'_1$ is equal to a  point  $y_o' \in \cV^{(2)} = g^{(1,2)}(\cV^{(1)})$.\par
Since $y_o \in \cV^{(1)}$ and $y_o' \in \cV^{(2)}$, there are exactly two points   $x_o \in \cU^{(1)} \cap \cE^n_{(p)}$ and $x_o' \in \cU^{(2)} \cap \cE^n_{(p)}$ such that 
$\varphi^{(p)}(x_o) = y_o$ and $\varphi^{(p)}(x'_o) = y'_o$. We may therefore consider  the unique $\cC^0$ curves $\h, \h': [0,1] \longrightarrow \cE^n_{(p)} \setminus Z_{F, \varphi^{(p)}}$ such that 
   \begin{itemize}
 \item[--] $\varphi^{(p)} \circ \h = \g$ and $\varphi^{(p)} \circ \h' = \g'$, 
 \item[--]   $\h_1 = x_o$ and $\h'_1 = x_o'$.  
  \end{itemize}
For any $t \in [0,1]$, consider the $F$-complete set $\{\h_t^{(1)} = \h_t, \h_t^{(2)}, \ldots, \h_t^{(m)} \}$ which contains $\h_t$. 
Then 
there exist $m$ neighbourhoods $\cU^{(j)}_t \subset \cE^n_{(p)}$, $1 \leq j \leq m$,  of the points $\h_t^{(j)}$ 
such that the restrictions
$$\varphi|_{\cU^{(j)}_t}: \cU^{(j)}_t \longrightarrow \cV^{(j)}_t = \varphi^{(p)}(\cU^{(j)}_t)\ ,\qquad F|_{\cU^{(j)}_t}: \cU^{(j)}_t \longrightarrow \cW_t\ ,\ \cW_t =  F(\cU^{(1)}_t)\ ,$$
are  biholomorphisms. Hence, also  the maps 
\beq \label{accacontbis}Êk^{(1,j)}_t =  \varphi^{(p)} \circ \left(F|_{\cU^{(j)}_t}\right)^{-1} \circ 
F|_{\cU^{(1)}_t} \circ \left(\varphi^{(p)}|_{\cU^{(1)}_t}\right)^{-1} : \cV^{(1)}_t \longrightarrow \cV^{(j)}_t\ ,\ \  2 \leq j \leq m\ ,\eeq
are biholomorphisms. Reordering  the elements in the $F$-complete set, we may always assume that 
 $$\cU^{(1)}_{t = 1} = \cU^{(1)} \cap \cE^n_{(p)}\ ,\qquad \cV^{(1)}_{t =1} = \cV^{(1)}\ ,\qquad k_{t = 1}^{(1,2)} = g^{(1,2)}|_{\cV^{(1)}}\ .$$
 By compactness and reorderings,  there exist $t_1$, \ldots, $t_N = 1$ $\in [0,1]$ such that $\g([0,1]) \subset \bigcup_{k = 1}^N \cV_{t_k}^{(1)}$ and 
$\cV^{(1)}_{t_j}Ê\cap  \cV^{(1)}_{t_{j-1}} \neq \emptyset$
for all $2 \leq j \leq N$. By the same arguments for \eqref{3.9}, we have  that 
$k^{(1,2)}_{t_j}|_{\cV^{(1)}_{t_j}Ê\cap \cV^{(1)}_{t_{j-1}}} = k^{(1, 2)}_{t_{j-1}}|_{\cV^{(1)}_{t_j}Ê\cap \cV^{(1)}_{t_{j-1}}}$ for all $2 \leq j \leq N$, so that
the map $k^{(1, 2)}_1$ extends to a holomorphic map  
$$k^{(1, 2)}:  \cV = \bigcup_{j = 1}^N \cV^{(1)}_{t_j}\longrightarrow  \cV' = \bigcup_{j = 1}^N \cV^{(2)}_{t_j}$$
between a neighbourhood $\cV$ of $\g([0,1])$ and a neighbourhood $\cV'$ of $k^{(1,2)}(\g([0,1]))$. 
Notice that, by construction, if  $\wt y$, $\wt y'$ are such that $\wt y' = k^{(1,2)}(\wt y)$,  they are  $F$-related.  
Since $k^{(1, 2)}|_{\cV^{(1)}_1} = g^{(1,2)}|_{\cV^{(1)}_1}$, by the Identity Principle,  
$k^{(1,2)}Ê= g^{(1,2)}|_{\cV}$ and  $y' = g^{(1,2)}(y) = k^{(1,2)}(y)$.  Therefore $y$, $y'$   are $F$-related, as we needed to prove. \par
\bigskip
{\it Case 2: $y$, $y' \in Z'$}. Let  $\{y_k\} \subset  B^n \setminus Z'$ be a sequence with 
 $\lim_{k \to \infty}Êy_k = y$. By continuity, the  sequence $y'_k =  g(y_k) $ converges to $y' = g(y)$. By the result in the previous case,  $y_k$ and $y'_k$ are $F$-related for any $k$ and there exists a sequence  $\{w_k\}  \subset  \O$ such that  $y_k, y'_k \in \varphi^{(p)}\left(F^{-1}(w_k)\right)$.
Since $\varphi^{(p)}$ and $F$ are proper, up to a subsequence,  we may assume   that  $\{w_k\}$ converges to a  point $w_o \in \O$. Using continuity, one can  check that this implies that $y, y'  \in \varphi^{(p)}(F^{-1}(w_o))$ and are therefore $F$-related. \par
\bigskip
Finally, the property that $\G$ is a subgroup follows from the fact that the composition of two elements $g^{(i,j)}$, $g^{(k, \ell)}$ $\in \G$ maps  the connected open set $\cV^{(1)}$ into one of the  $F$-related sets $\cV^{(r)}$. This can occur only if  $g^{(i,j)} \circ g^{(k, \ell)} |_{\cV^{(1)}}= g^{(1, r)}|_{\cV^{(1)}}$ for some $r$,  meaning that $g^{(i,j)} \circ g^{(k, \ell)} = g^{(1,r)} \in \G$. 
\end{pf}
\par
 \section{The Main Theorem}
\setcounter{equation}{0}
Consider now the proper holomorphic correspondence 
\beq\label{defpsi} \Psi = P_\G  \circ \varphi^{(p)} \circ F^{-1}: \Omega \longcorr B^n_{\G}\ , \eeq
where $P_\G$ and $B^n_{\G}$ are as defined in \S \ref{section2.5}. Theorem \ref{budino} is  
direct  consequence of the following:
\par
\begin{prop} The correspondence \eqref{defpsi} splits and each of its single-valued components  
$\Psi_i: \Omega \longrightarrow B^n_{\G}$, $1 \leq i \leq k$,  is a proper  holomorphic map such that 
\beq \label{numeretto}Ê \Psi_i \circ F=P_\G \circ \varphi^{(p)}\ .\eeq
\end{prop}
\begin{pf} By Lemma \ref{splittinglemma}, it suffices to show that the subset $S_\Psi \subset \O$ of the points $z$,  at which $\Psi$ does not split, 
is included in an analytic subvariety of dimension less than or equal to $n-2$. Let $\G_{\refl} \subset \G$ be the normal subgroup generated by the reflections in $\G$
and fix   some elements $h_1$, \ldots, $h_k $ in $\G \setminus \G_{\refl}$ such that $\G$ can be expressed as a disjoint union
 $$\G = \G_{\refl} \cup\G_{\refl}   h_1  \cup \ldots \cup \G_{\refl}   h_k\ . $$
 For convenience of notation, we set $h_0 = \Id_{B^n}$ so that $\G = \bigcup_{i = 0}^k \G_{	\refl} h_i$.\par
 \smallskip
 We first observe  that for any $g$, $g' \in \G_{\refl}$ and  $0 \leq i \neq  j \leq k$,  the element 
 $(g'  h_j)^{-1}Ê(g  h_i)$ is not in $\G_{\refl}$. In fact, since  $\G_{\refl}$ is normal,  if $\wt g =  h_j^{-1}Êg'{}^{-1}Êg  h_i$ is in $\G_{\refl}$, 
then 
 $$g'{}^{-1}Êg  h_i = h_j \wt g   = \wt g' h_j\ \text{for some}\ \wt g' \in \G_{\refl}\quad \Longrightarrow\quad \G_{\refl}Ê h_i \cap   \G_{\refl}Ê h_j \neq \emptyset\ ,$$
 contradicting the choice of the $h_m$'s. 
Due to this,  any fixed point set  $Fix((g'  h_j)^{-1}Ê(g  h_i))$ is an analytic variety of dimension less than or equal to $n-2$. \par
\smallskip
Let  $X$ be  the union of such fixed point sets, that is 
$$X = \bigcup_{\smallmatrix g, g' \in \G_\refl\\ 0 \leq i \neq j \leq k \endsmallmatrix}Fix((g'  h_j)^{-1}Ê(g  h_i))$$
 and note that 
$W = F(\varphi^{(p)}{}^{-1}(X))$ is  an analytic subvariety of $\O$ of dimension $\dim W \leq n-2$. In fact, $\wt X = \varphi^{(p)}{}^{-1}\left(X\right) \subset \cE^n_{(p)}$ is an analytic variety,   which is mapped  onto   $X$ and $W$ by the proper holomorphic maps $\varphi^{(p)}$ and $F$, respectively.  By the Proper Mapping Theorem, 
$$\dim W = \dim \wt X = \dim  X \leq  n-2\ .$$
\par
\medskip
Let  $w_o \in \O \setminus W$ and $z_o \in \varphi^{(p)} (F^{-1} (w_o))$. By construction, $z_o \notin X$. We claim that there exists a ball $B_{\varepsilon}(z_o) \subset B^n$, 
centred at $z_o$ and of radius $\varepsilon$, such that 
\beq \label{tina} g  h_i (B_{\varepsilon}(z_o)) \cap  g'  h_j (B_{\varepsilon}(z_o)) = \emptyset\eeq
for any $g$, $g' \in \G_{\refl}$ and $0 \leq i \neq j \leq k$. Suppose not.  Since $\G$ is finite,  there  exist $i \neq j$,  $g, g' \in \G_{\refl}$ and  two sequences $z_n$, $z'_n$  
such that 
$$z_o = \lim_{n \to \infty} z_n =  \lim_{n \to \infty} z'_n\qquad \text{and}\qquad  (g  h_i)(z_n) = (g'  h_j)(z_n') \ \text{for any}\ n\ .$$
 By continuity,   $z_o = ((g  h_i)^{-1}Ê (g'  h_j))(z_o)$,  i.e.  $z_o \in X$:  contradiction. \par
 In the following, we denote $ \cV_j = \bigcup_{g \in \G_{\refl}} g h_j(B_{\varepsilon}(z_o))$.  By \eqref{tina},  we have that  $\cV_i \cap \cV_j = \emptyset$ for any $0 \leq i \neq j \leq k$.\par
 \bigskip
We now consider an open  ball $B_{\d}(w_o) \subset \O$ with the following property: for any $w \in B_\d(w_o)$ there exists $z \in \varphi^{(p)}(F^{-1}(w))$ such that $z \in B_{\varepsilon}(z_o)$. The existence of such a ball can be checked as follows. Consider the $F$-complete set 
$F^{-1}(w_o)Ê= \{x^1_o, \ldots, x^N_o\}$ in  $ \cE^n_{(p)}$
and the corresponding $F$-complete set
$\varphi^{(p)}(F^{-1}(w_o)Ê)  = \{z_o^1, \ldots, z_o^{N'}\}Ê$ in  $B^n$. 
 Let $r$ sufficiently small so that $F^{-1}(B_r(w_o))$ has exactly $N$ connected components $U_1$, \ldots, $U_N$. Let $ V_i = \varphi^{(p)}(U_i)$ and assume that $z_o = z_o^1 \in V_1$. If there is no $B_\d(w_o)$ with the required property, there exists a sequence $\{w_\ell\} \subset  \O$, converging to $w_o$ such that 
 $$\varphi^{(p)}(F^{-1}(w_\ell)) \cap B_{\varepsilon}(z_o) = \emptyset\qquad \text{for any} \ \ell\ .$$
 Taking a suitable subsequence, we may assume that there exists a sequence $x_\ell \in U_1$ with $F(x_\ell) = w_\ell$ and $x_\ell$ converging to $x_o^1$. By construction,  the sequence 
 $\{\varphi^{(p)}(x_\ell) = z_\ell\}$ is in $  V_1$ and tends to  $\varphi^{(p)}(x_o^1) = z_o$. But this means that $z_\ell \in B_\varepsilon(z_o) \cap \varphi^{(p)}(F^{-1}(w_\ell))$
 for  all $\ell$'s sufficiently large and it contradicts our hypothesis. \par
 \medskip
 We now consider the  maps
  $$\psi_j: B_\d(w_o) \longrightarrow B^n_\G\ ,\ \psi_j(w) = P_\G(h_j(z))\  \text{for some} \ z \in \varphi^{(p)}(F^{-1}(w)) \cap \cV_0$$
with $0 \leq j \leq k$.   We claim that such maps  are well defined and single valued. In fact, if $z, z' \in \varphi^{(p)}(F^{-1}(w)) \cap \cV_0$, then,  by definition of $\cV_0$, 
  $$z = g (\wt z) \ ,\quad z' = g'(\wt z') \quad \text{for some}Ê\ \wt z, \wt z' \in B_{\varepsilon}(z_o)\ ,\ g, g' \in \G_{\refl} $$
  and, by Proposition \ref{proposition3.3}, $z' = h(z)$ for some $h \in \G$ and hence of  the form 
  $$ h =  g'' h_{i_o} \in \G_{\refl} h_{i_o} \quad \text{for some}Ê\ 0 \leq i_o \leq k\ .$$
 These two facts and the normality of $\G_\refl$ imply that 
$$g'(\wt z') = (g'' h_{i_o} g) (\wt z) =  (g'''  h_{i_o}) (\wt z)  \qquad \text{for some}\ g''' \in \G_\refl$$
 and hence that 
$$\wt z' = (\wh g h_{i_o}) (\wt z) \in \cV_{i_o}\qquad\text{with}\quad \wh g =  g'{}^{-1} g''' \in \G_\refl\ .$$
Since  $ \cV_0 \cap \cV_{i_o}= \emptyset$ for $i_o \neq 0$, we conclude that $h_{i_o} = h_0 =  \Id_{B^n}$ and  that $z' = g''(z) $. By normality of $\G_\refl$ and the properties of $P_\G$, it follows that 
$$P_\G(h_j (z')) = P_\G(h_j(g''(z))) = P_\G(g'''(h_j(z))) \overset{g''' \in \G_\refl}= P_\G(h_j(z))\ ,$$
 proving that $\psi_j$ is well defined and single valued.  Moreover, we have that 
\begin{lem} \label{lemmalemmetto}ÊEach map $\psi_j$ is  holomorphic.
\end{lem}
\begin{pf} Let us first show that the  $\psi_j$'s are continuous, i.e., that if 
 $w_\ell \in B_\d(w_o)$ is a sequence converging to $w \in B_\d(w_o)$, then $\lim_{\ell \to \infty} \psi_j(w_\ell) = \psi_j(w)$. Consider the $J$-complete  set $F^{-1}(w) = \{x^1, \ldots , x^N\} \subset \cE^n_{(p)}$. By construction of $B_\d(w_o)$, we may assume that 
  $z = \varphi^{(p)}(x^1)$ belongs to  $B_\varepsilon(z_o)$, so that $\psi_j(w) = P_\G(h_j(z))$. \par
Let $\overline{B_{r_i}(x^i)} \subset \cE^n_{(p)}$ be $N$ disjoint closed balls such that 
\beq  \label{4.2} F^{-1}(w) \cap \overline{B_{r_i}(x^{i})} = \{x^{i}\} \eeq
 and denote by 
$S \subset \O$ the compact set 
$S = \bigcup_{i = 1}^N F(\partial B_{r_i}(x^{(i)}))$.
Since $w \notin S$, there exists $B_{\d'}(w) \subset B_\d(w_o)$ such that $B_{\d'}(w) \cap S = \emptyset$. If we set 
\beq \label{4.3} R^{i} = F^{-1}\left(B_{\d'}(w) \right)\cap B_{r_i}(x^{i})\eeq
the  arguments of Prop. 15.1.6 in \cite{Ru} imply that the maps $F|_{R^{i}}: R^{i} \longrightarrow B_{\d'}(w)$ are proper and hence surjective. 
With no loss of generality, we may assume that $\{w_\ell\} \subset B_{\d'}(w)$ and we may consider a sequence   $\{x_\ell\} \subset R^{1}$ such that $F(x_\ell) = w_\ell$. Up to a subsequence,  $\{x_\ell\}$ converges to some $\wt x \in \overline{R^{1}}$. By \eqref{4.2},  \eqref{4.3} and continuity,   $F(\wt x) = w$ and $\wt x = x^1$. \par
Since    $\{z_\ell = \varphi^{(p)}(x_\ell)\} \subset B^n$  converges to $z = \varphi^{(p)}(x^1)$,  for  all $\ell$'s sufficiently large  $z_\ell$ is in $B_\varepsilon(z_o)$, so that
$\lim_{\ell \to \infty} \psi_j(w_\ell) = \lim_{\ell \to \infty}  P_\G(h_j(z_\ell)) =  P_\G(h_j(z)) = \psi_j(w)$, 
as claimed.
\par
\smallskip
We now prove that  $\psi_j$'s are holomorphic. In fact, for any $w \in B_\d(w_o) \setminus F(Z_F)$, there exist a neighbourhood $\cW$ of $w$ and 
neighbourhoods $\cU^1$, \ldots, $\cU^m$ of the pre-images $x^1$, \ldots, $x^m$ of $w$, such that 
$F|_{\cU^i}: \cU^i \longrightarrow F(\cU^i) = \cW$
are biholomorphisms. For any $z \in \varphi^{(p)}(F^{-1}(w)) \cap B_{\varepsilon}(z_o)$,  there exists $1 \leq j_o \leq m$ such that   
$$z = \varphi^{(p)}(F|_{\cU^{j_o}}^{-1}(w))\ .$$ 
Taking  $\cW$ sufficiently small, we may suppose that for any $w' \in \cW$
$$z' = \varphi^{(p)}(F|_{\cU^{j_o}}^{-1}(w')) \in B_{\varepsilon}(z_o) \quad \Longrightarrow \quad \psi_j(w') = P_\G\circ h_j\circ \varphi^{(p)}\circ F|_{\cU^{j_o}}^{-1}(w') \ ,$$
proving that $\psi_j|_{\cW}$ is holomorphic.  This implies that 
$\psi_j$ is holomorphic in $B_\d(w_o) \setminus F(Z_F)$ and, by continuity and known facts on holomorphic extensions (\cite{Ru},  Cor. of Thm. 4.4.7), it  is holomorphic on   $B_\d(w_o)$. 
\end{pf}
\par
By construction,  for any $w$ $\in B_\d(w_o)$ we have that $\Psi(w) =$ $(\psi_0(w),$ $\psi_1(w), \ldots, \psi_k(w))$. By Lemma \ref{lemmalemmetto}, the $\psi_j$'s are holomorphic, meaning that  $\Psi$ splits at $w_o$. Since $w_o$ is an arbitrary point of $\O \setminus W$ and $\dim W \leq n-2$, by Lemma \ref{splittinglemma} we have that $\Psi$ splits. The equality  \eqref{numeretto}  is a direct consequence of the definition of $\Psi$. 
\end{pf}


    \end{document}